\documentclass{amsart}
\usepackage{amssymb, array, verbatim, amscd}
\usepackage{amsmath, amsfonts,amsthm}

\newtheorem{theorem}{Theorem}[section]

\newtheorem{definition}[theorem]{Definition}
\newtheorem{example}[theorem]{Example}

\newcommand{\Z}{{\mathbb Z}}
\newcommand{\R}{{\mathbb R}}
\begin{document}

\title{Probabilistic Gene Regulatory Networks, isomorphisms of Markov Chains }

\author{Mar\'{\i}a Alicia Avi\~n\'o }
\address{ Department of
Mathematic-Physics,
       University of Puerto Rico,
       Cayey, PR 00736}
\thanks{This research was supported by the National Institute of Health,
 PROGRAM SCORE, 2004-08, 546112, University of Puerto Rico-Rio Piedras Campus, IDEA Network of Biomedical
  Research Excellence, and the Laboratory Gauss University of Puerto Rico Research.
   I want to thank Professor E. Dougherty for his useful suggestions, and Professor O. Moreno for  his support
   during the last four years. }
\email{
mavino@cayey.upr.edu}
 \date{\today}

 \begin{abstract} In this paper we study homomorphisms of  Probabilistic Regulatory Gene Networks(PRN) introduced in \cite{A1}. The model PRN is a natural generalization of the
  Probabilistic Boolean Networks (PBN), introduced by I. Shmulevich, E. Dougherty, and W. Zhang in \cite{SDZ},  that has been using to describe genetic
 networks and has therapeutic applications, see \cite{SGHDZ}.
In this paper, our main objectives are to apply the  concept of homomorphism and
$\epsilon$-homomorphism of probabilistic regulatory networks to the dynamic of the networks. The meaning of $\epsilon$ is that these homomorphic networks have similar distributions and the distance between the   distributions is upper bounded by $\epsilon$.Additionally, we prove that  the class of
PRN together with the homomorphisms form a category with products and coproducts.
 Projections are special homomorphisms, and they always induce
invariant subnetworks that contain all the cycles and steady  states in the network. Here, it is
proved that the $\epsilon$-homomorphism for $0<\epsilon<1$ produce simultaneous Markov Chains in both networks, that permit to introduce the concept of  $\epsilon$-isomorphism of Markov Chains, and similar networks.
\end{abstract}
\maketitle \section*{Introduction}
Genes can be understanding in their complexity behavior using models according with their discrete or continuous action. Developing computational tools permits describe gene functions and understand
the mechanism of regulation \cite{48,49}. \emph{This understanding will have a significant impact on the development of techniques for drugs testing and therapeutic intervention  for treating human diseases}\cite{D2005,50,SGHDZ}.

We focus our attention in the discrete structure of  genetic regulatory networks, instead of, its dual moving continuo-discrete.
Probabilistic Gene Regulatory Network(PRgN)  is  a natural generalizations of the model Probabilistic Boolean Network (PBN), introduced by I. Shmulevich, E. Dougherty, and W. Zhang in \cite{SDZ}. The mathematical background of the model PgRN, is introduced here, for simplicity we work with functions defined over  a set $X$ to itself, with probabilities assigned to these functions. $X$ is a set of states of genes, for example $X=\{0,1\}^n$, if our network is a Boolean network. Working in this way, we can observe the dynamic of the network indeed focus our attention in the description of functions. The set $X$ can be a subset of $\{0,1\}^n$, and we can extend some classical ideas to regulatory network, such as invariant subnetworks, automorphisms group, etc.  In particular if $X$  is a vector space over a finite field, the functions are lineal functions, then we can use linear algebra to describe the state space. Mapping are important in the study of networks, because they permit to recognize subnetworks, in particular determine when two networks are similar or equivalent. Special mappings are homomorphisms and $\epsilon$-homomorphisms, we use both to describe subnetworks and similar networks. An homomorphism transform a network to another in such a way the discrete structure giving by the first network can  lives in part of the other one, or these two networks are very similar but no equals, in particular in the probabilistic way. An $\epsilon$-homomorphism is the same but  with the condition that the probability distributions of the networks are close, and we use a preestablishes $0<\epsilon<1 $  as a distance between the probabilities. For concept of homomorphism of discrete dynamical systems see \cite{ABM,ABM1,LP}.
\section{Preliminaries concepts:Finite dynamical systems, probabilistic Boolean networks and Probabilistic Regulatory Networks}
Two finite dynamical systems $(X,f)$ and $(Y,g)$ are isomorphicas(or equivalents) if there exists a bijection $\phi: X\rightarrow
Y$ such that $\phi\circ f=g\circ \phi$, ( or $ f=\phi^{-1}\circ
g\circ \phi)$. If $\phi$ is not a bijection map then $\phi$ is an
homomorphism.

If $Y\subset X$ is such that $f(Y)\subset Y$ then $(Y,f|_y)$ is a
sub-FDS of $(X,f)$, where $f|_Y$ is the map restricted to Y. There exists naturally
  an injective morphism from $Y$ to $X$ called inclusion and denoted by
  $\iota$. The state space of a FDS $(X,f)$, is a digraph whit vertices the
  set $X$, and with an arrow from $u$ to $v$ if $f(u)=v$.

For example, the  FDSs
$X=(\{0,1\}^2,f_1(x,y)=(xy,y))$, and
$Y=(\{0,1\}^2,f_2(x,y)=(x,(x+1)y))$ are isomorphics, because their
state spaces are isomorphics.
\[\begin{array}{cccccc}
 (1,0)&\rightarrow ^{f_1}  & (0,0)&  & (0,1)&\\
    & &\circlearrowright & &\circlearrowleft ^{f_1} & \curvearrowright^{f_1} \\
 & & & & &(1,1) \cr
\end{array} \]
\[\begin{array}{cccccc}
 (1,1)&\rightarrow ^{f_2}  & (1,0)&  & (0,0)&\\
    & &  \circlearrowright & &\circlearrowleft ^{f_2} & \curvearrowright^{f_2} \\
 & & & & &(0,1) \cr
\end{array} \]
In fact, the isomorphism $\phi:\{0,1\}^2\rightarrow \{0,1\}^2$ is  the
bijection $\phi(1,0)=(1,1)$, $\phi(0,0)=(1,0)$, $\phi(0,1)=(0,0)$,
and $\phi(1,1)=(0,1)$. The following is an example of homomorphism
(inclusion) with $Z=\{\{(0,0),(1,0)\},f_1\}$.
\[\begin{array}{ccccccccc}
(1,0)&&  & (1,0)& & &  & (0,1)&\\
  \downarrow f_1&&\hookrightarrow ^\iota & \downarrow f_1 & &  & &\circlearrowleft ^{f_1} & \curvearrowright^{f_1} \\
 (0,0)&\circlearrowleft&&(0,0)& \circlearrowleft& & & &(1,1) \cr
\end{array} \]
A Probabilistic Boolean Network $\mathcal{A}=(V,F,C)$ is defined by
the following sort (type) of objects \cite{SDZ,SDKZ}:
 a set of nodes (genes) $V=\{x_1,\ldots ,x_n\}$, $x_i\in
\{0,1\}$, for all $i$; a family  $F=\{F_1, F_2,\ldots,F_n\}$ of
ordered sets $
F_i=\{f_1^{(i)},f_2^{(i)},\ldots,f_{\ell(i)}^{(i)}\}$ of Boolean
functions $f_j^{(i)}:\{0,1\}^n\rightarrow\{0,1\}$, for all  $j$
called predictors;  and a list $C=(C_1, \ldots ,C_n)$,
$C_i=\{c_1^{(i)},\ldots ,c_{\ell(i)}^{(i)}\}$, of selection
probabilities. The selection probability that the function
$f_j^{(i)}$ is used for the vertex $i$ is
$c_j^{(i)}=Pr\{f^{(i)}=f_j^{(i)}\}$. The dynamic of the PBN is
given  by a vector of functions $\mathbf{f}
_k=(f_{k_1}^{(1)},f_{k_2}^{(2)},\ldots,f_{k_n}^{(n)})$ for
$1\le{k_i}\le{ l(i)},$ and $f_{k_i}^{(i)}\in{F_i}$, where
$k=[k_1,\ldots , k_n]$, $1\leq k_i\leq \l (i).$ The map
 $\mathbf{f}_k: {\{0,1\}}^n\rightarrow{\{0,1\}} ^n$
acts as a transition function. Each variable $x_i\in\{0,1\}^n$
represents the state of the vertex $i$.  All functions are updated
synchronously. At every time step, one of the functions is
selected randomly from the set $F_i$ according to a predefined
probability distribution. The selection probability that the
transition function  $\mathbf{f}
_k=(f_{k_1}^{(1)},f_{k_2}^{(2)},\ldots,f_{k_n}^{(n)})$ is used to
 go  from the state $u\in \{0,1\}$ to another state $\mathbf{f}
_k(u)=v\in \{0,1\}^n$ is given by
 \[c_{\mathbf{f}_k}=\prod_{i=1}^n c_{k_i}^{(i)}.
\]
The dynamical transition structure of a PBN can be described by a
Markov chain with fixed transition probabilities. There are two
digraphs structures associated with a PBN: the low-level digraph
$\Gamma$, consisting of genes functions essentiality relations;
and the high-level digraph which consists of the states of the
system and the transitions between states. The  matrix $T$
associated to the high level digraph formed by placing $p(u,v)$ in
row $u$ and column $v$, where $u,\ v\  \in \{0,1\}^n$ is called
the transition probability matrix or chain matrix,
$p(u,v)=\sum_{\mathbf{f_k| f_k}(u)=v}c_{\mathbf{f}_k}$.

\subsection{Probabilistic Regulatory Gene Networks}
 A Probabilistic Gene Regulatory Network (PRN) is
a triple $\mathcal{X}=(X,F,C)$ where $X$ is a finite set and
$F=\{f_1, \ldots , f_n\}$ is a set of functions from $X$ into
itself, with a list $C=(c_1, \ldots ,c_n)$ of selection
probabilities, where $c_i=p(f_i)$, \cite{A1,A}
 We associate with each PRN  a weighted digraph, whose vertices are the elements of $X$, and
if $u,v\in X$, there is an arrow going from $u$ to $v$ for each
function $f_i$ such that $f_i(u)=v$, and the probability $c_i$ is
assigned to this arrow. This weighted digraph will be called the state
space of $\mathcal{X}$. In this paper, we use the notation PRN for one or more networks.\\
\begin{example}
If $X=\{0,1\}^2$,
$F=\{f_1(x,y)=(x,y), f_2(x,y)=(x,0),$\\$f_3(x,y)=(1,y),f_4(x,y)=(1,0)\}$; and $C=\{.46,.21,.22,.11\}$,   the state space  of
$\overline{\mathcal{X}}=(X,F,C)$ is  the following:\\
\[\begin{array}{c}
\overset{.67}{\circlearrowright}(0,0)\leftarrow ^{.21} (0,1)\circlearrowleft^{.46}\\
^{.33}\downarrow \hbox{   }\hbox{  }\swarrow_{ .11}\downarrow^{.22}\\
\overset{1}{\circlearrowright}(1,0)\overset{.32}{\longleftarrow} (1,1)\circlearrowleft^{.68} \cr
\end{array} T=\left[\begin{array}{cccc} .67&0&.33&0\\
.21&.46&.11&.22\\
0&0&1&0\\
0&0&.32&.68\cr
\end{array}\right]\]
\end{example}
\subsection{Homomorphisms and $\epsilon$-homomorphisms of PRN}
If $C$ is a set of selection probabilities we denote by $\chi$ the
characteristic function over $C$. That is $\chi:C\cup\{0\}
\rightarrow \{0,1\}$ such that $\chi(c)=1$, if $c\ne 0$ and
$\chi(0)=0$. Let
$\mathcal{X}_1=(X_1,F=(f_i)_{i=1}^n ,C)$ and
$\mathcal{X}_2=(X_2,G=(g_j)_{j=1}^m, D)$ be two PRN.\\
\begin{definition}[ Homomorphisms of PRN] A  map $\phi:X_1\rightarrow X_2$ is  an
\textbf{homomorphism} from $\mathcal{X}_1$ to $
\mathcal{X}_2$, if for all $f_i$ there exists a $g_j$, such that
for all $u$, $v$ in $\mathcal{X}_1$,
\[\hbox{(1) } \phi \circ f_i=g_j\circ \phi; \hbox{ and  } \hbox{(2) }\chi(d_{g_j}(\phi(u),\phi(v)))\geq
\chi(c_{f_i}(u,v)).\]
\[\begin{array}{c}
X_1\overset{f_i}{\longrightarrow} X_1\\
\phi \downarrow  \hspace{.4in} \downarrow \phi \\
X_2\overset{g_j}{\longrightarrow} X_2\cr
\end{array}\]
(3) Condition for $\epsilon$-Homomorphism: The
distributions of probabilities following the homomorphism are
enough close. An $\epsilon$- homomorphism is an
homomorphism  that satisfies the condition, for all $i$, $j$, $max|p(u_i,u_j)-p(\phi(u_i),\phi(u_j))|\le \epsilon$, where $\epsilon >0$ is  a real number that we previously  determine for the  applications.

If $\phi:X_1\rightarrow X_2$ is a bijective map, and for all $f_i$, $g_j$, $u$, and $v$ in $\mathcal{X}_1$;  $d_{g_j}(\phi(u),\phi(v))=c_{f_i}(u,v)$, then $\phi$  is an isomorphism.
\end{definition}
  \begin{example}
  \end{example}
If    $\overline{\mathcal{X}}=(X;F;C)$ is the PRN in Example $1$, and  $\mathcal{X}_1=(X;F'=\{f_1, f_2,f_3\};C'=\{.47,.28,.25\})$ is a new PRN over the same set $X$ with  different probabilities and only three functions.
\[\begin{array}{c}
\mathcal{X}_1\\
^{.75}\circlearrowright(0,0)\leftarrow ^{.28} (0,1)\overset{.47}{\circlearrowleft}\\
^{.25}\downarrow \hspace{.4in}\hbox{  }\downarrow^{.25}\\
\circlearrowright^{1}(1,0)\overset{.28}{\longleftarrow }(1,1)\circlearrowleft^{.72} \cr
\end{array}  \overset{\phi}{\hookrightarrow}    \begin{array}{c}
\overline{\mathcal{X}}\\
^{.67}\circlearrowright(0,0)\leftarrow ^{.21} (0,1)\overset{.46}{\circlearrowleft}\\
^{.33}\downarrow \hbox{   }\hbox{  }\swarrow_{ .11}\downarrow^{.22}\\
\circlearrowright^{1}(1,0)\overset{.32}{\longleftarrow} (1,1)\circlearrowleft^{.68} \cr
\end{array}\]
\[T_1=\left[\begin{array}{cccc} .75&0&.25&0\\
.28&.47&0&.25\\
0&0&1&0\\
0&0&.28&.72\cr
\end{array}\right],\overline{T}=\left[\begin{array}{cccc} .67&0&.33&0\\
.21&.46&.11&.22\\
0&0&1&0\\
0&0&.32&.68\cr
\end{array}\right]\]
The homomorphism $\phi : \mathcal{X}_1\rightarrow \overline{\mathcal{X}}$ is a bijective map, $\phi (x)=x$, over the set of states, but an inclusion over the set of arrows, because the arrow going from  $(0,1)$ to $(1,1)$ in $\overline{\mathcal{X}}$ doesn't appear in  $\mathcal{X}_1$.
 The first condition for homomorphism is obvious. The condition (2) holds, because the inclusion of arrows. The two transition matrices are connected by this inclusion, since  if the place $ij$ in the first matrix $\ne 0$ then this place is $\ne 0$  in the second network too. The two PRN are not isomorphics because the probabilities are not equals. In order to determine $\epsilon$ for the homomorphism, we use the transition matrices. In this example $\epsilon =.11$.\\
\[T_1-\overline{T}=\left[\begin{array}{cccc} .08&0&-.08&0\\
.07&.01&-.11&.03\\
0&0&0&0\\
0&0&-.04&.04\cr
\end{array}\right]\]

If the homomorphism is a bijective map like here, the transition matrices $T_1$ and $T_2$ have the same order, and   $\sum_{i=1}^n (T_1-T_2)_{ij}=0$, for $j=\bar{1,n}$\\
\section{Applications to Markov Chains, $\epsilon$-Similar Networks }

 Two PRN  are $\epsilon$-similar if  there exists a bijective homomorphism  $\phi$ between them, such that $\phi ^{-1}$ is also an homomorphism. Observe that $\phi$ and $\phi^{-1}$ have the same $\epsilon$.

\begin{example}
\end{example}The following networks are $.005$-homomorphics  with the injective homeomorphism
$\phi:\{0,1\}^2\rightarrow \{0,1\}^3$, giving by $\phi(x,y)=(x,y,1)$. We can see in this example, that the network $\mathcal{X}_1$ is an $\epsilon$-subnetwork of $\mathcal{X}_2$,  the special subnetwork $\phi (\mathcal{X}_1)$ is an invariant subnetwork of $\mathcal{X}_2 $. The networks, and the transition matrices are the following:
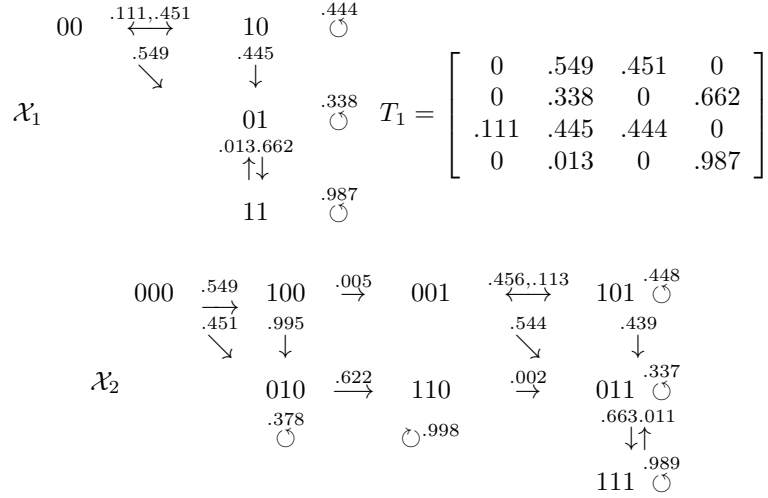
\[\mathcal{X}_1 \begin{array}{cccc}
00&\overset{.111, .451}{\longleftrightarrow}&10&\overset{.444}{\circlearrowleft }\\
&\overset{.549}{\searrow }&\overset{.445}{\downarrow }&\\
&&{01}&\overset{.338}{\circlearrowleft }\\
&&\overset{.013 .662}{\uparrow \downarrow }&\\
&& 11&\overset{.987}{\circlearrowleft }\cr
\end{array}\  T_1=\left[\begin{array}{cccc} 0&.549&.451&0\\
0&.338&0&.662\\
.111&.445&.444&0\\
0&.013&0&.987\cr
\end{array}
\right]\]
\[\mathcal{X}_2\begin{array}{ccccccc}
000&\underrightarrow{^{.549}}&100&\overset{.005}{\rightarrow }&001&\overset{.456,.113}{\longleftrightarrow}&101\overset{.448}{\circlearrowleft }\\
\hspace{.2in}&\overset{.451}{\searrow }&\overset{.995}{\downarrow }&&&\overset{.544}{\searrow }&\overset{.439}{\downarrow }\\
&&010&\overset{.622}{\longrightarrow }&110&\overset{.002}{\rightarrow }&011\overset{.337}{\circlearrowleft }\\
&&\overset{.378}{\circlearrowleft}&&\circlearrowright ^{.998}&&\overset{.663.011}{\downarrow \uparrow }\\
&&&&&& 111\overset{.989}{\circlearrowleft }\cr
\end{array}\]
Ordering the elements in the following way
\[\{(0,0,0),(0,1,0),(1,0,0),(1,1,0),(0,0,1),(0,1,1),(1,0,1),(1,1,1)\}\]
\[T_2=\left[\begin{array}{cccccccc}
0&.451&.549&0&0&0&0&0\\
0&.378&0&.622&0&0&0&0\\
0&.995&0&0&.005&0&0&0\\
0&0&0&.998&0&.002&0&0\\
0&0&0&0&0&.544&.456&0\\
0&0&0&0&0&.337&0&.663\\
0&0&0&0&.113&.448&.439&0\\
0&0&0&0&0&.011&0&.989\cr
\end{array}.
\right]\]

Where, $ T_\phi=\left[\begin{array}{cccc} 0&.544&.456&0\\
0&.337&0&.663\\
.113&.448&.439&0\\
0&.011&0&.989\cr
\end{array}
\right]$. Observe that,
\[T_1-T_{\phi}=\left[\begin{array}{cccc} 0&.005&-.005&0\\
0&.001&0&-.001\\
-.002&-.003&.005&0\\
0&.002&0&-.002\cr
\end{array}
\right]\]
As a consequence, we obtain $max|(T_1)_{ij}-(T_\phi)_{ij}|\leq .005$, then the networks are $.006$-homeomorphiscs. The steady state of $T_1$ is $\pi_1=(0,.01926,0,.98074)$, and the steady
state of $T_\phi$ is $\pi_\phi=(0,.01632,0,.98368)$. We can see that $|\pi_1-\pi_\phi|=max_i|\pi_1(i)-\pi_\phi(i)|<.004$. Meanwhile for $T_2$ we
have, $\pi_2=(0,0,0,.01632,0,0,0,.98368)$.
We can observe that with the new order in the elements of ${\Z_2}^3$, the transition matrix  $T_\phi$ is an invariant  submatrix of $T_2$, that is $T_2=\left[\begin{array}{cc}
T_{11} & T_{12} \\
0& T_\phi\cr
\end{array}
\right]$.
In the above example,  $\mathcal{X}_1$ and $\phi(\mathcal{X}_1)$ are $\epsilon$-similar.

To introduce the following theorem we use the transition matrices in this example, and we can observe the following calculation
\[T_1^2-{T_\phi}^2=\left[\begin{array}{cccc}
-.001467&-.00136&.00006&.00277\\
0&.00 199&0&-.00 199\\
-.00232&-.00019&.00295&-.00243\\
0&.002639& 0&-.00 263
\cr
\end{array}
\right],\]
therefore $max|(T_1^2)_{ij}-(T_\phi^2)_{ij}|\leq .003$.
\[T_1^3-{T_\phi}^3=\left[\begin{array}{cccc}
-.000394&-.00044&.00011&.00073\\
0&.002525&0&-.00253\\
-.000161&.00156&.00213&-.00353\\
0&.002843&0&-.002843\cr
\end{array}
\right],\]
and $max|(T_1^3)_{ij}-(T_\phi^3)_{ij}|\leq .004$.
\begin{theorem}\label{teo}
If $\phi:\mathcal{X}_1\rightarrow\mathcal{X}_2$ is an $\epsilon$-homomorphism, then the transition matrices $T_1$ and $T_\phi$ satisfy the condition
\[max|({T_1}^n)_{ij}-({T_\phi}^n)_{ij}|\leq \epsilon\]
 for all possible $i$ and $j$, and all $n>1$. If the homomorphism is injective and $\epsilon< 1 $, the steady state of $T_1$ and the steady state of $T_\phi$ are close, that is satisfy $|\pi_1-\pi_\phi|=max_i|\pi_1(i)-\pi_\phi(i)|\leq\epsilon$.
 \end{theorem}
\begin{proof}
Since  $({T_1}^n)_{ij}=p(u_i,f^n(u_i))$ when $f^n(u_i)=u_j$, therefore we have the following:
\[|p(u,f^2(u))-p(\phi(u),\phi(f^2(u))|=\]
\[|p(u,f(u))p(f(u),f^2(u))-p(\phi(u),\phi(f(u)))p(\phi(f(u)),\phi(f^2(u)))|\leq\]
\[|p(f(u),f^2(u))||p(u,f(u))-p(\phi(u),\phi(f(u)))|+\]
\[|p(\phi(u),\phi(f(u)))||p(f(u),f^2(u))-p(\phi(f(u)),\phi(f^2(u)))|\leq\epsilon\]
Using this property, mathematical induction  over $n$, and the cases giving for the  possible values of the probabilities $p(u,v)$, we can conclude that our  aim holds.
\end{proof}
An $\epsilon$-homomorphism between two PRN determines a correspondence between the Markov Chains of these two  networks. The following  new definition is connected in some way with the idea of bisimulation for two TDMC (Time Discrete Markov Chains) studied and introduced in \cite{BS}. Two TDMC are equivalent by a bisimulation if they simulate the same  problem with two different set of probabilities. Here, we introduce the concept of two similar TDMC
\begin{definition} The two TDMC of the same size $n\times n$: $\{T_1,\ T_1^2, \ T_1^3,\  \ldots\}$, and $\{T_2,\ T_2^2, \ T_2^3, \ \ldots\}$ are $\epsilon$-similar  or $\epsilon$-isomorphic if
\begin{itemize}
\item[(1)] there exists an $\epsilon \in \R$ enough small, such that $T_1^m-T_2^m=(t_{ij})_{n\times n}$ satisfies that $|t_{ij}|<\epsilon$, and $\sum_{i=1}^n t_{ij}=0$, for all $m$,
\item[(2)] $\chi(T_1^m)_{ij}=\chi(T_2^m)_{ij}$,  for all $m$, where $\chi$ is the characteristic function.
\end{itemize}
That is,  simulating the networks by these two TDMC are $\epsilon$-similar.
\end{definition}
In the above example, the TDMC generated by $T$ and $T_\phi$ are $.005$-similar, and the networks simulated by them are $.005$-similar.
 \section{Construction of Probabilistic Regulatory Networks, examples}
 \medskip
\par\noindent
\textbf{Sum of two PRN}

 Let $\mathcal{X}_1=(X_1,F=(f_i)_{i=1}^n ,C)$ and
$\mathcal{X}_2=(X_2,G=(g_j)_{j=1}^m, D)$ be two PRN. The sum
$\mathcal{X}_1\oplus\mathcal{X}_2=(X_1\dot{\cup} X_2, F\vee G,C
\vee D)$ is a  PRN where
\begin{itemize}
\item  [(1)]$X_1\dot{\cup} X_2$ is the disjoint union of $X_1$ and
$X_2$.
 \item [(2)]the
function $h_{ij}=(f_i\vee g_j)$ is defined by $h_{ij}(x)=f_i(x)$
if $x\in X_1$ and $h_{ij}(x)=g_j(x)$ if $x\in X_2$.

\item [(3)]the probability $p(h_{ij})=c_i \vee d_j$, that is
$p(h_{ij})=c_i$ if $h_{ij}=f_i$ or $p(h_{ij})=d_j$ if
$h_{ij}=g_j$.
\end{itemize}

If $T_1$ and $T_2$ are the transition matrices of $\mathcal{X}_1$ and $\mathcal{X}_2$ respectively,
Then $T=\left( \begin{array}{cc}
T_1 &0\\
0& T_2\cr
\end{array}\right)$ is the transition matrix of $\mathcal{X}_1\oplus\mathcal{X}_2$.
\begin{example}
\end{example} An example of sum is the PRN obtained by summing
the same PRN twice, $\mathcal{X}\oplus \mathcal{X}$. To make the
disjoint union, we  subindicate  $X$ with $0$ for the first $X$
and with $1$ for the second $X$. That is, the new set is
\[X_0\dot{\cup} X_1=\{(0,0,0),(0,1,0),(0,1,0),(1,1,0)\}\]
\[\cup \{(0,0,1),(0,1,1),(0,1,1),(1,1,1)\}.\]
The digraph is:
\[\begin{array}{ccccc}
   \curvearrowright^{.6} & & \curvearrowright ^{.4} & &\curvearrowright ^{1} \\
 (1,1,0)&\rightarrow ^{.4} & (1,0,0)& \rightarrow ^{.6} & (0,0,0)\\
    & &  & &  \curvearrowright ^{1}\\
 & & & & (0,1,0) \cr
\end{array}\]
\[\begin{array}{ccccc}
   \curvearrowright^{.6} & & \curvearrowright ^{.4} & &\curvearrowright ^{1} \\
 (1,1,1)&\rightarrow ^{.4} & (1,0,1)& \rightarrow ^{.6} & (0,0,1)\\
    & &  & &  \curvearrowright ^{1}\\
 & & & & (0,1,1) \cr
\end{array}  \]

This is a way to construct a PRN over $\{0,1\}^n$ using either one
or two  PRN over $\{0,1\}^{n-1}$, since $2^{n-1}+2^{n-1}=2^n$.
\medskip
\par\noindent
\textbf{Superposition}

It is clear that a PRN is the
\textbf{superposition} of several Finite dynamical Systems (FDS)\cite{H} over the same set $X$ with
probabilities assigned to each FDS. Since each functions defined
over a finite field can be wrote as a polynomial function, we will
use this notation for functions over a finite field, \cite{AGM}. If
$X=\{0,1\}=\Z_2$, the finite field of two elements,  all the FDSs
over $X$ have one of the following state space, where $f_1(x)=x;\
f_2(x)=1;\ f_3(x)=0;\ f_4(x)=x+1$, $\forall x\in X$:

\[\begin{array}{c}
L_1\\
   0\circlearrowleft\\
 \\
1\circlearrowleft\cr
\end{array}
\begin{array}{c}
L_2\\
0\\
 \downarrow\\
1\circlearrowleft\cr
\end{array}
\begin{array}{c}
L_3\\
0\circlearrowleft\\
 \uparrow\\
 1\cr
\end{array}  \begin{array}{c}
L_4\\
0\\
\uparrow \downarrow\\
1\cr

\end{array}  \]
If $p_i$ denotes de probability assigned to $L_i$, and $T_i$ denotes its transition matrix,
 then the set of all PRN is described as follows.
\[\left\{ (X,F,C)|T=\sum _{i=1}^4 p_iT_i=\left(\begin{array}{cc}
   p_1+p_3&p_2+p_4\\
p_3+p_4& p_1+p_2\cr
\end{array} \right) \sum _{i=1}^4 p_i=1
\right\}\]
We denote by $L_1L_2$ the superposition of $L_1$ and $L_2$, and similarly  $L_1L_3$ is the superposition of $L_1$ and $L_3$. The following state spaces are the superposition of two FDS with two elements:
 \[\begin{array}{c}
 L_1L_2\\
   0\circlearrowleft^{p_1}\\
 \downarrow ^{p_2}\\
1\circlearrowleft^{1}\cr
\end{array},  \  \   \begin{array}{c}
L_1L_3\\
0\circlearrowleft^{1}\\
 \uparrow^{p_3}\\
 1 \circlearrowleft^{p_1}\cr
\end{array} , \  \begin{array}{c}
L_1L_4\\
0\circlearrowleft^{p_1}\\
 \overset{p_4}{\uparrow \downarrow}\\
 1 \circlearrowleft^{p_1}\cr
\end{array} ,  \
\begin{array}{c}
L_2L_3\\
   0\circlearrowleft^{p_3}\\
 \downarrow ^{p_2}\uparrow^{p_3}\\
1\circlearrowleft^{p_3}\cr
\end{array}  \  \   \begin{array}{c}
L_2L_4\\
\hspace{-.3in}0\\
 \downarrow^{p_4} \uparrow  ^{p_2+p_4}\\
 1 \circlearrowleft^{p_3}\cr
\end{array} , \
\begin{array}{c}
L_3L_4\\
0\circlearrowleft^{p_3}\\
 \overset{p_4}{\downarrow}\uparrow^ {p_4+p_3}\\
\hspace{-.2in} 1 \cr
\end{array},  \]

For example, with transition matrices \[T_{12}=T_1+T_2=\left(\begin{array}{cc}
   p_1&p_2\\
0& 1\cr
\end{array} \right) \  T_{13}=T_1+T_3=\left(\begin{array}{cc}
   1&0\\
p_3& p_1\cr
\end{array} \right)\]
\medskip
\par\noindent
\textbf{Product of two PRN}

Let
$\mathcal{X}_1=(X_1,F=(f_i)_{i=1}^n ,C)$ and
$\mathcal{X}_2=(X_2,G=(g_j)_{j=1}^m, D)$ be two PRN. The product
$\mathcal{X}_1\times \mathcal{X}_2=(X_1\times X_2, F \times G,C
\wedge D)$ is a  PRN where
\begin{itemize}
\item  [(1)]$X_1\times X_2$ is the cartesian product  of $X_1$ and
$X_2$.
 \item [(2)]the
function $h_{ij}=(f_i,g_j)$ is defined by
\[h_{ij}(x_1,x_2)=(f_i(x_1),g_j(x_2))\]  for $x_1\in X_1$, and
$x_2\in X_2$.

\item [(3)]the probability $p(h_{ij})$ is a function of $c_i$ and $d_j$, for example $p(h_{ij})=\frac{c_i+d_j}{2}$.
\end{itemize}
\begin{example}
\end{example}
The product $L_1L_2\times L_1L_3$ is the PRN with four states $\{(0,0),(0,1),$
$(1,0),(1,1)\}$ and
four functions \[f_{11}(x,y)=(x,y), f_{13}(x,y)=(x,0),\]
 \[ f_{21}(x,y)=(1,y),\  f_{23}(x,y)=(1,0).\]
The state space is the following:
\[\begin{array}{c}
L_1L_2\times L_1L_3\\
^{p_{11}+p_{13}}\circlearrowright(0,0)\leftarrow ^{p_{13}} (0,1)\circlearrowleft^{p_{11}}\\
^{p_{23}+p_{21}}\downarrow \hbox{   }\hbox{  }\swarrow_{ p_{23}}\downarrow^{p_{21}}\\
\hbox{      }\hbox{      }\hbox{ }\hbox{  }\hbox{   }\hspace{.3in}    \hbox{     } \  \hbox{    }  \circlearrowright^{1}(1,0)\longleftarrow (1,1)\circlearrowleft^{p_{11}+p_{21}} \\
\hspace{.3in} ^{p_{13}+p_{23}} \cr
\end{array}\]
The transition matrix is the following
\[T=\left(\begin{array}{cccc}
 p_{11}+p_{13}  &0&p_{23}+p_{21}&0\\
p_{13}& p_{11}&p_{23}&p_{21}\\
0&0&1&0\\
0&0&p_{13}+p_{23}& p_{11}+p_{21}\cr
\end{array} \right)\]
\subsection{Linear Probabilistic Regulatory Networks }

 A linear PRN is a
\textbf{\emph{superposition}} of linear FDS. A linear FDS is a
pair $(X,f)$ where $f$ is a linear function, and $X$ is a vector
space over a finite field. So, a linear PRN is a triple $(X,(f_i)_{i=1}^m, C)$, where $X$ is
a finite vector space, the functions $f_i:X\rightarrow X$ are linear functions, and
 $C=\{c_i=p(f_i)\}$. The set $X$ has cardinality a power of a prime number and each linear function
 is determined by its characteristic polynomial and the companion matrix.

  If $X=\Z_3=\{ 0,1,2\}$ is the field  of
integer modulo $3$, then the linear functions are:  $f_1(x)=x$,
$f_2(x)=2x$, and $f_3(x)=0$ for all $x\in \Z_3$. So, the  linear
PRN are the following:
\[\hbox{   }\begin{array}{c}
\{ f_1,f_2\}\\
  \circlearrowleft ^{1}\\
   \mathbf{0} \\
   \\
  \circlearrowright^{p_1} \mathbf{1}\rightleftarrows ^{p_2}\mathbf{2}\circlearrowleft^{p_1}\cr
\end{array} \hbox{   }
\hbox{    }\begin{array}{c}
\{f_1,f_3\}\\
\circlearrowleft ^{1}\\
   \mathbf{0} \\
\nearrow^{p_3}\nwarrow  \\
\mathbf{1} \circlearrowleft _{p_1}
 \hbox{    }\circlearrowright_{p_1} \mathbf{2}
\cr
\end{array} \
\begin{array}{c}
\{f_2,f_3\}\\
\circlearrowleft ^{1}\\
   \mathbf{0} \\
\nearrow^{p_3}\nwarrow  \\
\mathbf{1} \rightleftarrows _{p_2}  \mathbf{2}
\cr
\end{array}
\hbox{    }\begin{array}{c}
\{f_1,f_2,f_3\}\\
\circlearrowleft ^{1}\\
   \mathbf{0} \\
\nearrow^{p_3}\nwarrow  \\
\circlearrowleft _{p_1}\mathbf{1} \rightleftarrows _{p_2}  \mathbf{2}\circlearrowleft _{p_1}
\cr
\end{array}\]
If $X=\Z_2\times \Z_2$ is the vector space with $4$ elements over
the field $\Z_2$, then there are $4$ linear FDS not isomorphics. In fact,
using matrix, the possible characteristics polynomials $p_f(\lambda)$ are:
$\lambda ^2,\ \lambda ^2+\lambda \ \lambda ^2 +1,\ \lambda
^2+\lambda +1$. The companion matrices of these linear functions are:
\[A_1=\left(\begin{array}{cc}
0&0\\
0 &0\cr
\end{array}\right) \hbox{   }A_2=\left(\begin{array}{cc}
0&0\\
0 &1\cr
\end{array}\right) \  A_3=\left(\begin{array}{cc}
1&0\\
0 &1\cr
\end{array}\right)\hbox{   }A_4=\left(\begin{array}{cc}
0&1\\
1 &1\cr
\end{array}\right)\]
Then the FDS associated to this matrices are:
\[\begin{array}{ccc}
A_1&&\\
\circlearrowright(0,0)&\leftarrow& (1,0)\\
\uparrow&\nwarrow&\\
(0,1)&&(1,1)\cr
\end{array}
\hbox{   }\begin{array}{ccc}
A_2&&\\
\circlearrowright(0,0)&\leftarrow& (1,0)\\
&&\\
\circlearrowright(0,1)&\leftarrow&(1,1)\cr
\end{array}\hbox{   }\]
\[\begin{array}{ccc}
A_3&&\\
\circlearrowright(0,0)&& \circlearrowright(1,0)\\
&&\\
\circlearrowright(0,1)&&\circlearrowright(1,1)\cr
\end{array}\hbox{    }\begin{array}{ccc}
A_4&&\\
\circlearrowright(0,0)&& (1,0)\\
&\swarrow&\uparrow\\
(0,1)&\rightarrow&(1,1)\cr
\end{array}\hbox{   }\]
The linear PRN  with two functions are the following:
\[\begin{array}{c}
A_1,A_2\\
{\overset{1}{\circlearrowright}}(0,0){\overset{1}{\leftarrow}(1,0)} \\
^{p_1}\uparrow \nwarrow ^{p_1}\\
{\overset{p_2}{\circlearrowright}}(0,1){\overset{p_2}{\leftarrow}}(1,1)\cr
\end{array}
\begin{array}{c}
A_1,A_3\\
{\overset{1}{\circlearrowright}}(0,0){\overset{p_3}{\leftarrow}}(1,0){\overset{p_3}{\circlearrowleft}} \\
^{p_1}\uparrow\nwarrow^ {p_1}\\
^{p_3} \circlearrowright\hbox{  }(0,1) \   (1,1)\circlearrowleft ^{p_3}\cr
\end{array} \
 \begin{array}{c}
A_1,A_4\\
\circlearrowright^1(0,0){\overset{p_1}{\leftarrow}}(1,0)\\
^{p_1}\uparrow \hspace{.1in}  \overset{\nwarrow}{\swarrow}\hspace{.1in} \uparrow^{p_4}\\
\  (0,1)  {\overset{p_4}{   \leftarrow}}(1,1)
\cr
\end{array}\]
\[\begin{array}{c}
A_2,A_3\\
\circlearrowright^{1}(0,0)\leftarrow^{p_2} (1,0)\circlearrowleft ^{p_3}\\
\\
\circlearrowright^{1}(0,1)\leftarrow^{p_2}(1,1)\circlearrowleft ^{p_3}\cr
\end{array}
\begin{array}{c}
A_2,A_4\\
\circlearrowright^1(0,0)\leftarrow ^{p_2} (1,0)\\
\hspace{.2in}^{p_4}\swarrow \uparrow ^{p_4}\\
\circlearrowright^{p_2}(0,1)\leftarrow ^{p_4}(1,1)\cr
\end{array}\hbox{    }
\begin{array}{c}
A_3,A_4\\
\circlearrowright^1 (0,0)\  \   (1,0)\circlearrowleft ^{p_3}\\
\    ^{p_4}\swarrow\uparrow ^{p_4} \\
^{p_3}\circlearrowright(0,1)\rightarrow^{p_4}(1,1)\circlearrowleft^{p_3}
\cr
\end{array}\hbox{   }
\]
\begin{example} Looking the state spaces of the $A_i,A_j$ for $i\ne j$ PRN, we can conclude that two of them  are not isomorphic. Working with the transition matrices, we can calculate when there exists  an homomorphism between two of them. For example: the set $Hom(A_1A_2;A_1A_3)=\emptyset$
\end{example}
\begin{example}
\end{example}
 If $X_1$ and $X_2$ are two PRNs then the product $X_1\times X_2$ can be projected over each component. In the example of product  of two PRN, we use the PRN $L_1L_2$ and $L_1L_3$ to construct $L_1L_2\times L_1L_3$. The functions $\phi_1(x,y)=x$ and $\phi_2(x,y)=y$, the usual projections, are $\epsilon$-homomorphism. In fact

 (1)$\phi_i:\Z_2\times \Z_2\rightarrow \Z_2$, and  $\phi_1\circ f_{11}=f_1\circ \phi_1$,and $\phi_1\circ f_{21}=f_2\circ \phi_1$.

 (2) the second condition is satisfies too, since  several arrows compact to one arrow in $\Z_2$, for example
 $\phi_1(0,1)=\phi_1(0,0)=0$, then the arrows compact in the following way
 \[\begin{array}{ccc}
 \{^{p_{11}+p_{13}}\circlearrowright(0,0)\leftarrow ^{p_{13}} (0,1)\circlearrowleft^{p_{11}}\}& \twoheadrightarrow & 0\circlearrowleft ^{p_1}\cr
 \end{array}\]
 Similarly
 \[\begin{array}{ccc}
 \{\circlearrowright ^1(1,0)\leftarrow ^{p_{13}+p_{23}} (1,1)\circlearrowleft^{p_{11}+p_{21}}\}& \twoheadrightarrow & 1\circlearrowleft {1}\cr
 \end{array}\]
 \[\left\{\begin{array}{c}
 (0,0)\hspace{.6in} (0,1)\\
 \downarrow^{p_{23}+p_{21} } \hspace {.1in} \swarrow_{p_{23}} \downarrow ^{p_{21}}\\
 (1,0) \hspace {.6in}(1,1)\cr
 \end{array}\right\}\twoheadrightarrow
  \begin{array}{c}
  0\hspace{.4in}\\
  \hspace{.1in}\downarrow_{p_2}\\
  1\hspace{.4in}\cr
  \end{array}\]

  (3) for the $\epsilon$-condition, the probabilities need to satisfy:
  \[|p_{11}-p_1|,\ |p_{13}-p_1|,\ |(p_{11}+p_{13})-p_1|<\epsilon\]
  \[ |(p_{13}+p_{23})-1|,\ |(p_{11}+p_{21})-1|<\epsilon\]
  \[|p_{21}-p_2|,\ |p_{23}-p_2|,\ |(p_{21}+p_{23})-p_2|<\epsilon\]
  since $p_1+p_2=1$, and $p_1+p_3=1$,  for $p_{ij}=p_ip_j$ for example, we have  $|(p_{11}+p_{13})-p_1|=|(p_{21}+p_{23})-p_2|=0$, $|p_{11}-p_1|=|p_{23}-p_2|=|p_1p_2|$, $|p_{13}-p_1|=|p_1^2|$,  $|p_{21}-p_2|=|p_2^2|$, $|(p_{13}+p_{23})-1|=|p_1|$, and $|(p_{11}+p_{21})-1|=|p_2|=|p_3|$. So, max$(p_1,p_2,p_3)=\epsilon$.

  We have the same conditions for $\phi_2$. It is clear that always the projections are homomorphisms, but  the applications of  the $\epsilon$-condition depend on the particular case that we are working. It is clear that there are several ways to define inclusions $\iota :L_1L_3\rightarrow L_1L_3\times L_1L_2$.

  The following digraph shows the projections in general.
  \[\begin{array}{c}
  X_1\times X_2\\
  \swarrow_{\pi_1}\hspace{.1in} _{\pi_2}\searrow\\
  X_1\hspace{.6in}X_2\cr
  \end{array}\]
  The functions $\pi_i:X_1\times X_2\rightarrow X_i$,  are defined are follows: $\pi_i(x_1,x_2)=x_i$, $i=1,2$. It is easy to check that the others two conditions are satisfied as in the above example.

\section{Invariant Subnetworks and Projections}
A subnetwork  $Y\subseteq X$  of $\mathcal{X}=(X,F,C)$ is an
\textbf{invariant subnetwork or a sub-PRN} of $\mathcal{X}$ if
 $f_i(u)\in Y$ for all $u\in Y$, and  $f_i\in F$. Sub-PRNs
 are sections of a PRN, where there aren't arrows going out.
  The complete network $X$, and any cyclic state with
  probability 1, are sub-PRNs. An invariant subnetwork is irreducible if doesn't have a
proper invariant subnetwork.
 \emph{An endomorphism is a projection if $\pi^2=\pi$. }\\
\begin{theorem}
If there  exists a projection from $\mathcal{X}$ to a subnetwork $\mathcal{Y}$ then $\mathcal{Y}$ is  an invariant subnetwork of $\mathcal{X}$.
\end{theorem}
\begin{proof}
Suppose that there exists a projection $\pi:X\rightarrow Y$. If  $y\in Y$,  by definition of projection $\pi(y)=y$,  and $f_i(\pi(y))=\pi(g_j(y))$. Therefore all  arrows in the subnetwork $Y$ are going inside $Y$, and the network is invariant.
\end{proof}

\begin{example}
\end{example}
The PRN $\overline{\mathcal{X}}$ has two invariant subnetworks with projections \\
 $\pi_1(x,y)=(x,0)$ and $\pi_2(x,y)=(1,y)$.
\[\overline{S}_1\begin{array}{ccc}
&&\hspace{.4in}\\
\overset{.67}{\circlearrowright}0 &&\overset{.67}{\circlearrowright}(0,0)\\
^{.33}\downarrow &\cong &^{.33}\downarrow   \\
\overset{1}{\circlearrowright} 1& &\overset{1}{\circlearrowright}(1,0)\cr
\end{array} S_1\overset{\pi_1}{\longleftarrow}\begin{array}{c}
 \hspace{.3in}\overset{.67}{\circlearrowright}(0,0)\leftarrow ^{.21} (0,1)\overset{.46}{\circlearrowleft}\\
\hspace{.3in}^{.33}\downarrow \hbox{   }\overset{.11}{\swarrow}\downarrow^{.22}\\
\hspace{.3in}\circlearrowright^{1}(1,0)\overset{.32}{\longleftarrow }(1,1)\overset{.68}{\circlearrowleft} \\
\hspace{-1.3in}\overline{\mathcal{X}}\\
\hspace{-1.3in}\overset{\pi_2}{\downarrow}\\
S_2\circlearrowright^{1}(1,0)\overset{.32}{\longleftarrow} (1,1)\overset{.68}{\circlearrowleft}\hspace{-.3in}\\
\cong\hspace{-.3in} \\
\overline{S}_2\circlearrowright^{1}0\overset{.32}{\longleftarrow} 1\overset{.68}{\circlearrowleft}\hspace{-.3in}\cr
\end{array} \]
Checking the probabilities for $\pi_1$ and $\pi_2$, we have
$\epsilon_1=.68$; and  $\epsilon_2=.67$.
We can observe that
$\overline{\mathcal{X}}\cong \overline{S}_1\times \overline{S}_2$.
\begin{example}
\end{example} The subnetwork $\mathcal{X}_1=(\{(x,y,1)\}, F, C)$ is an invariant subnetwork of $\mathcal{X}=(\{0,1\}^3,F, C) $.
\[\mathcal{X}\begin{array}{ccccccc}
&&&&&\mathcal{X}_1&\\
000&\underrightarrow{^{.549}}&100&\overset{.005}{\rightarrow }&\textbf{001}&\overset{\textbf{.113,.456}}{\longleftrightarrow}&\textbf{101}\overset{\textbf{.448}}{\circlearrowleft} \\
\hspace{.2in}&\overset{.451}{\searrow }&\overset{.995}{\downarrow }&&&\overset{\textbf{.544}}{\searrow }&\overset{\textbf{.439}}{\downarrow} \\
&&010&\overset{.622}{\longrightarrow }&110&\overset{.002}{\rightarrow }&\textbf{011}\overset{\textbf{.337}}{\circlearrowleft }\\
&&\overset{.378}{\circlearrowleft}&&\circlearrowright ^{.998}&&\overset{\textbf{.663.011}}{\downarrow \uparrow }\\
\overset{\pi}{\downarrow}&&&&&& \textbf{111}\overset{\textbf{.989}}{\circlearrowleft }\cr
\end{array}\]
\[\mathcal{X}_1 \overset{\rho}{\cong} \overline{\mathcal{X}_1}\begin{array}{cccc}
00&\overset{.113, .456}{\longleftrightarrow}&10&\overset{.448}{\circlearrowleft }\\
&\overset{.544}{\searrow }&\overset{.439}{\downarrow }&\\
&&{01}&\overset{.337}{\circlearrowleft }\\
&&\overset{.011 .663}{\uparrow \downarrow }&\\
&& 11&\overset{.989}{\circlearrowleft }\cr
\end{array}\]
Ordering the elements in the following way $\{(0,0,0),(0,1,0),\\(1,0,0),(1,1,0),(0,0,1),(0,1,1),(1,0,1),(1,1,1)\}$, the matrix\\
 $ T_{X_1}=\left[\begin{array}{cccc} 0&.544&.456&0\\
0&.337&0&.663\\
.113&.448&.439&0\\
0&.011&0&.989\cr
\end{array}
\right]$ is an invariant part of the transition matrix
 $T_{\mathcal{X}}=\left[\begin{array}{cc}
T_{11}& T_{12}\\
0&T_{X_1}\cr
\end{array}\right].$\\
Using the projection $\pi:\mathcal{X}\rightarrow \mathcal{X}_1$,
$\pi(x,y,z)=(x,y,1)$; and the isomorphism $\rho (x,y,1)=(x,y)$,
the network $\mathcal{X}$ is projected over the network
$\overline{\mathcal{X}_1}$. Checking the arrows the projection
$\overline{\pi}$ is a  $.5$-homomorphism.
\begin{example}
 The following PBN appears  in \cite{ID}, and has three sub-PBNs.
 \end{example}
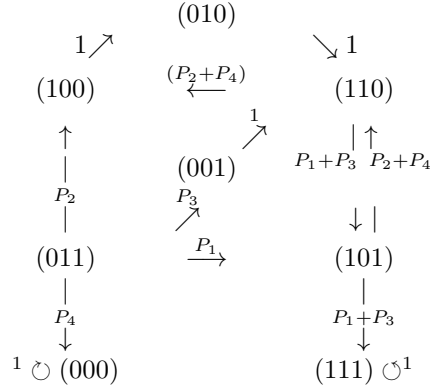
\[ \begin{array}{cccc}
&(010)&&\\
\hspace{.3in}1{\nearrow}&&\hspace{-.3in}\searrow 1\hspace{.6in}\\
(100)&\overset{(P_2 +P_4)}{\longleftarrow} &(110)\\
\uparrow&\hspace {.5in}\overset{1}{\nearrow} & | \uparrow  &\\
|&(001)&^{P_1+P_3} \hspace{.05in} ^{P_2+P_4}& \\
\overset{P_2}{|}&\hspace{-.2in}\overset{P_3}{\nearrow} &{\downarrow}\ {|}&\\
(011)&\overset{P_1}{\longrightarrow} &(101)&\\
|&&|&\\
\overset{P_4}{\downarrow}&&\overset{P_1+P_3}{\downarrow}&\\
^ 1\circlearrowright(000)&&(111)\circlearrowleft ^1\cr
\end{array}\]
 $X_1=\{(000)\}$, $X_2=\{(111)\}$, and
 $X_3=\{(100),(010),(110),(101),(111)\}$ are sub networks.
With adequate order the transition matrices $T_X$ and $T_{X_3}$ are:
\[T_{X_3}=\left[\begin{array}{ccccc} 0&0&.5&.5&0\\
1&0&0&0&0\\
0&1&0&0&0\\
.5&0&0&0&.5\\
0&0&0&0&1\cr
\end{array}
\right]\hbox{ and }\pi_{X_3}=(0,0,0,0,1)\]
meanwhile $T_X=\left[\begin{array}{cc}
T_{11} & T_{12} \\
0& T_{X_3}\cr
\end{array}\right]$, and $\pi_X=(.25,0,0,0,0,0,0,.75)$.

\section{ The category of Probabilistic Regulatory Networks, and mathematical background}
\begin{theorem}
  If $\phi_1:\mathcal{X}_1\rightarrow \mathcal{X}_2$ is an
$\epsilon_1$-homomorphism, and
$\phi_2:\mathcal{X}_2\rightarrow \mathcal{X}_3$ is another $\epsilon_2$-homomorphism. Then $\phi=\phi_2\circ \phi_1:\mathcal{X}_1\rightarrow \mathcal{X}_3$ is an $\epsilon$-homomorphism. Therefore the Probabilistic Regulatory Networks with the homomorphisms of PRN form the category \textbf{PRN}.
\end{theorem}
\begin{proof}
The Probabilistic Regulatory Networks with the PRN homomorphisms  is a
category if: the composition is an homomorphism, and satisfy the
associativity law; and there exists an identity homomorphism for
each PRN.

  (1) Let $\phi_1:\mathcal{X}_1\rightarrow \mathcal{X}_2$ be an $\epsilon_1$-homomorphism, and let
$\phi_2:\mathcal{X}_2\rightarrow \mathcal{X}_3$ be an $\epsilon_2$-homomorphism.
If $q_t$, $g_k$ and $f_j$ are
functions in each PRN, and such that $\phi_1 \circ f_j=g_k \circ
\phi_1$ and $\phi_2 \circ g_k=q_t \circ \phi_2$, then we
will  prove that: $\phi \circ f_j=q_t \circ \phi.$ In fact,
 \[(\phi_2 \circ \phi_1) \circ f_j=\phi_2 \circ (\phi_1 \circ f_j)=\]
 \[\phi_2 \circ (g_k
 \circ \phi_1)=(\phi_2 \circ  g_k)\circ \phi_1=\]
 \[(q_t \circ \phi_2)\circ \phi_1=q_t \circ (\phi_2\circ \phi_1).\]
\noindent
 (2)We want to prove that
\[
\chi(t_k(\phi(u),\phi(v)))\ge{\chi(c_i(u,v))}.
\]
Suppose that
${\chi(c_i(u,v))}=1$.
Then, since $\phi_1$ is an homomorphism of PRN, we have that
\[\chi(d_j (\phi_1(u),\phi_1(v)))\ge {\chi (c_i(u,v))%
}\]
 which is 1. Since
$\phi_2$ is an homomorphism of PRN, we obtain that
\[\chi(t_k(\phi(u),
\phi(v)))=\chi(t_k(\phi_2(\phi_1(u)), \phi_2(\phi_1(v))))\]

\[ \ge \chi(c_j(\phi_1(u),(\phi_1(v))=1.\]
Therefore we obtain that
\[
\chi (t_k(\phi_2(\phi_1(u)), \phi_2(\phi_1(v))))=1.
\]

Then the composition of two PRN-homomorphisms  is an homomorphism.

(3) To verify the third condition for $\epsilon$-homomorphism, we do the following.
If $p(\phi(u_1),\phi(u_2))>1$, with $u_1,\ u_2\in X_1$, then we need to prove that there exists an $\epsilon$ such that \[|p(u_1,u_2)-p(\phi(u_1),\phi(u_2))|<\epsilon.\] In fact:

\[|p(u_1,u_2)-p(\phi(u_1),\phi(u_2))|=|p(u_1,u_2)-p(\phi_1(u_1),\phi_1(u_2))+ \] \[p(\phi_1(u_1),\phi_1(u_2))-p(\phi_2(\phi_1(u_1)),\phi_2(\phi_1(u_2)))|\]
\[< |p(u_1,u_2)-p(\phi_1(u_1),\phi_1(u_2))|+\]
\[|p(\phi_1(u_1),\phi_1(u_2))-p(\phi_2(\phi_1(u_1)),\phi_2(\phi_1(u_2)))|\leq \epsilon_1+\epsilon _2\]
\[|p(u_1,u_2)-p(\phi(u_1),\phi(u_2))|< \epsilon_1+\epsilon _2\]
because $\phi_1$ and $\phi_2$ are $\epsilon$-homomorphisms.

The associativity and identity laws are easily checked, therefore
our claim holds, and \textbf{PRN} is a category.
\end{proof}
It is clear that, the PRN with the homomorphism between them form a category that we will denote $\mathcal{PRN}$. The category \textbf{PRN} is a  subcategory of $\mathcal{PRN}$,  since an homomorphism is not always an homomorphism for some $\epsilon \in \R$ enough small. But, if we don't include the condition for $\epsilon$ to be enough small, the two categories are the same, because always an homomorphism is an $\epsilon$-homomorphism for some $ \epsilon \in \R$.
\begin{theorem}\label{product}
Let $\mathcal{X}_1\times \mathcal{X}_2=(X_1\times X_2,H,E)$ be a product of PRN $\mathcal{X}_1=(X_1,F,C)$ and $\mathcal{X}_2=(X_2,G,D)$. If $\delta_i:X \rightarrow X_i$ are two PRN-homomorphisms, then there exists an homomorphism $\delta: X\rightarrow X_1\times X_2$, such that $\phi_i\circ \delta=\delta_i$ for $i=1,2$. That is, the following diagram commutes
\[\begin{array}{c}
  \hspace{.05in} X_1\times X_2 \hspace{.05in}\\
 \overset{\phi_1}{ \swarrow}\hspace{.1in}\overset{\delta}{\uparrow}\hspace{.1in}\overset{\phi_2}{\searrow}\\
  X_1\overset{\delta_1}{\longleftarrow} X \overset{\delta_2}{\longrightarrow} X_2\cr
  \end{array}\]
  This homomorphism is unique.
\end{theorem}
\begin{proof}
The function $\delta:X\rightarrow X_1\times X_2$ is defined as follows $\delta(x)=(\delta_1(x),\delta_2(x))$,  $x\in X$. $\delta$ is an homomorphism, in fact:

(1) Let $\mathcal{X}=(X,L,P)$ be a PRN. Since $\delta_1$ and $\delta_2$ are homomorphism, for all function $l_t\in L$ there exist two functions $f_i\in F$ and $g_j\in G$, such that $\delta_1\circ l_t=f_i\circ \delta_1$, and $\delta_2\circ l_t=g_j\circ \delta_2$. Then for the function $l_t$ there exists the function $(f_i,g_j)$ that satisfies $\delta \circ l_t=(f_i,g_j)\circ \delta$.
\[(\delta \circ l_t)(x)=\delta(l_t(x))=(\delta_1(l_t(x)),\delta_2(l_t(x)))=\]
\[(f_i(\delta_1(x)),g_j(\delta_2(x)))=((f_i,g_j)\circ\delta)(x)\]

(2) In order to prove $\chi(e_{ij}(\delta(x),\delta(x')))\geq \chi(p_{l_t}(x,x'))$,  suppose  $\chi(p_{l_t}(x,x'))=1$. Then $l_t(x)=x'$, and $\delta(x')=\delta(l_t(x))=(f_i,g_j)(\delta(x))$ by  part (1). Therefore $\chi(e_{ij}(\delta(x),\delta(x')))=1$, and our claim holds.

It is easy to check that $\phi_i\circ \delta=\delta_i$, in fact \[\phi_1 (\delta(x))=\phi_1(\delta_1(x),\delta_2(x))=\delta_1(x),\] for all $x\in X$.
\end{proof}
If $\delta_i$, $i=1,2$,  are $\epsilon _i$-homomorphism then \[max|p(x,x')-p(\phi_1(\delta(x)),\phi_1(\delta(x')))|\leq \epsilon_1.\] But
\[|p(x,x')-p(\delta(x),\delta(x'))+p(\delta(x),\delta(x'))-
p(\phi_1(\delta(x)),\phi_1(\delta(x')))|\leq \]
\[|p(x,x')-p(\delta(x),\delta(x'))|+\]
\[|p(\delta(x),\delta(x'))-
p(\phi_1(\delta(x)),\phi_1(\delta(x')))|\leq \epsilon_1.\]
Therefore
\[|p(x,x')-p(\delta(x),\delta(x'))|\leq \epsilon_1-|p(\delta(x),\delta(x'))-
p(\phi_1(\delta(x)),\phi_1(\delta(x')))|\]
\[|p(x,x')-p(\delta(x),\delta(x'))|\leq \epsilon_1-\overline{\epsilon}_1.\]
Therefore $\delta$ is an $\epsilon$-homomorphism. So, the theorem holds for $\epsilon$-homomorphism.

It is an immediate consequence the following result, also is true for $\epsilon$-homomorphisms.
\begin{theorem}
Let $\mathcal{X}_1\oplus \mathcal{X}_2=(X_1\times X_2,H,E)$ be a product of PRN $\mathcal{X}_1=(X_1,F,C)$ and $\mathcal{X}_2=(X_2,G,D)$. If $\gamma_i:X_i \rightarrow X$ are two PRN-homomorphisms, then there exists an homomorphism $\gamma: X_1\oplus X_2\rightarrow X $, such that $\gamma \circ \iota_i=\gamma_i$ for $i=1,2$. That is, the following diagram commutes
\[\begin{array}{c}
  \hspace{.05in} X_1\oplus X_2 \hspace{.05in}\\
 \overset{\iota_1}{ \nearrow}\hspace{.1in}\overset{\gamma}{\downarrow}\hspace{.1in}\overset{\iota_2}{\nwarrow}\\
  X_1\overset{\gamma_1}{\longrightarrow} X \overset{\gamma_2}{\longleftarrow} X_2\cr
  \end{array}\]
  This homomorphism is unique.
\end{theorem}

\begin{theorem} All  reducible PRN is either a product of its non trivial sub-PRN  or a subnetwork of this product.
\end{theorem}
\begin{proof} It is trivial by definition of Product and sub-PRN.
\end{proof}
\section{Conclusions}
The intersection, and the union of two sub-PRN  is a sub-PRN, therefore the class of sub-PRN of a particular PRN is a lattice. Reduction mappings described in \cite{ID} and defined for PBN using the influence of a gene, for example $x_n$, on the predictor function $f_j^{(i)}$ to determine the selected predictor, can be extended to PRN. In order to extend this procedure to more than boolean functions, we use the polynomial description of genetic functions given in \cite{AGM}, the partial derivative is the usual in calculus  and all the concepts  in \cite{ID} can be using for PRN. Similarly our definition of projection, the reduction mappings are $\epsilon$-homomorphisms, and we can use  for genes with more than two quantization, since this extension is not a trivial work we  develop the theory and methods  in \cite{A}.

\end{document}